\documentclass[12pt]{article}
\usepackage{amssymb}
\def\bsk{\bigskip}
\def\ds{_{(s)}}
\def\dv{_{(v)}}
\def\du{_{(u)}}
\def\dr{_{(r)}}
\def\bp{{\bar p}}
\def\k{\kappa}

\newcommand{\eq}{\begin{equation}}
\newcommand{\en}{\end{equation}}
\newcommand{\eqs}{\begin{eqnarray*}}
\newcommand{\ens}{\end{eqnarray*}}
\newcommand{\eqa}{\begin{eqnarray}}
\newcommand{\ena}{\end{eqnarray}}

\newcommand{\ex}{\mathbb E}

\def \proof{\noindent{\it Proof.\ }}

\newtheorem{theorem}{\large Theorem}[section]

\newtheorem{lemma}[theorem]{\large  Lemma}

\def\numberlikeadb{\global\def\theequation{\thesection.\arabic{equation}}}
\numberlikeadb

\def\ignore#1{}

\def\Z{{\mathbb Z}}

\def\nn{{\cal N}}
\def\Eq{\ =\ }
\def\Le{\ \le\ }
\def\Ref#1{(\ref{#1})}

\def\sji{\sum_{j\ge1}}
\def\dtv{d_{{\rm TV}}}
\def\Le{\ \le\ }
\def\Ge{\ \ge\ }
\def\law{{\cal L}}

\def\nin{\noindent}

\def\Bi{{\rm Bi\,}}
\def\MN{{\rm MN\,}}
\def\TP{{\rm TP\,}}
\def\Po{{\rm Po\,}}
\def\pr{{\bf P}}
\def\half{{\textstyle{{1\over2}}}}

\def\b{\beta}
\def\e{\varepsilon}
\def\z{\zeta}
\def\slo{\sum_{l\ge0}}

\def\law{{\cal L}}

\def\ep{\hfill$\Box$}
\def\nn{{\cal N}}
\def\var{{\rm Var\,}}
\def\cov{{\rm Cov\,}}
\def\s{\sigma}
\def\Blb{\Bigl\{}
\def\Brb{\Bigr\}}
\def\Blm{\left|}
\def\Brm{\right|}

\def\ul{^{(l)}}

\def\non{\nonumber}
\def\a{\alpha}
\def\r{\rho}

\def\d{\delta}

\def\etal{{\it et~al.}}
\def\m{\mu}
\def\n{\nu}

\def\l{\lambda}
\def\t{\tau}
\def\Bl{\Bigl(}
\def\Br{\Bigr)}

\def\sjjn{\sum_{j\ge j_n}}
\def\ssjn{\sum_{j\ge j_n}}
\def\ssj{\sum_{s\ge j}}
\def\stj{\sum_{t\ge j}}
\def\dloc{d_{{\rm loc}}}
\def\dtv{d_{{\rm TV}}}
\def\fdd{\|f''\|}
\def\wjnl{U_j^{(n-l)}}
\def\wjn{U_j^{(n)}}
\def\mjnl{M_{j\cdot}^{(n-l)}}
\def\mjn{M_{j\cdot}^{(n)}}
\def\zjl{Z_{j\cdot}^{(l)}}
\def\giv{\,|\,}

\def\ssjnj{\sum_{{s\ge j_n \atop s\ne j}}}
\def\ssjn{\sum_{s\ge j_n}}
\def\stjn{\sum_{t\ge j_n}}

\def\um{^{(m)}}
\def\tyl{\t^{-1}y_j(l)}

\def\tssq{\t^{-1}\sjjn\slo q_j(l)}
\def\tssqy{\t^{-1}\sjjn\slo  q_j(l) y_j(l)}
\def\ssq{\sjjn\slo q_j(l)}
\def\stt{\s^2\t^{-3}}
\def\sqy{\slo  q_j(l) y_j(l)}
\def\sq{\slo  q_j(l)}
\def\slq{\sq}

\def\zjsl{Z_{js}\ul}
\def\mjsnl{M_{js}^{(n-l)}}
\def\mjsn{M_{js}^{(n)}}
\def\hm{{\hat \m}}

\begin{document}

\title{Univariate approximations in the infinite occupancy scheme\footnote
{A.D. Barbour gratefully acknowledges financial support from Schweizerischer
Nationalfonds Projekt Nr.~20-117625/1.}
}
\author{A. D. Barbour\thanks{Angewandte Mathematik,
Winterthurerstrasse~190, CH--8057 Z\"urich, Switzerland:\hfil\break
{\tt a.d.barbour@math.uzh.ch}}
\\
{University of Z{\"u}rich}} % and Utrecht University}}
\date{\empty}

\maketitle

\begin{abstract}
\noindent
In the classical occupancy
scheme with infinitely many boxes, $n$ balls are thrown independently into
boxes $1,2,\ldots$, with probabilities~$p_j$, $j\ge1$.
We establish approximations to the distributions of the summary
statistics~$K_n$, the number of
occupied boxes, and~$K_{n,r}$, the number of boxes containing exactly~$r$ balls,
within the family of translated Poisson distributions.
These are shown to be of ideal order as $n\to\infty$,
with respect both to total variation distance and to the approximation
of point probabilities.  The proof is probabilistic, making use of
a translated Poisson approximation theorem of R\"ollin~(2005).
\end{abstract}

\hglue2.6cm
local limit approximation


\vskip1.2cm

\section{Introduction}
\setcounter{equation}{0}
In the classical occupancy
scheme with infinitely many boxes,  $n$ balls are thrown independently into
boxes $1,2,\ldots$, with probability $p_j$ of hitting box $j$, $j\ge1$, where
$p_1\geq p_2\geq\ldots>0$ and $\sum_{j=1}^\infty p_j=1$.
The summary statistics~$K_n$, the number of
occupied boxes, and~$K_{n,r}$, the number of boxes containing exactly~$r$ balls,
have been widely studied.  Central limit theorems were established by
Karlin~(1967), under a regular variation condition, and Dutko~(1989) showed
that~$K_n$ is asymptotically normal,
assuming only the necessary condition that its variance
tends to infinity with~$n$. A full discussion of this and many more
aspects of the problem can be found in Gnedin~\etal\ (2007); see also
Barbour \& Gnedin~(2009), in which multivariate approximation of the~$K_{n,r}$
is treated.

As regards the accuracy of the central limit approximation, Hwang \& Janson~(2008)
show that the point probabilities $\pr[K_n=t]$ are uniformly approximated by
the point probabilities of the integer discretization of the normal distribution
$\nn(\m_n,\s_n^2)$, where $\m_n := \ex K_n$ and $\s^2_n := \var K_n$.
The accuracy of their approximation is  of order
$O(1/\s^2_n)$, provided only that $\s^2_n \to \infty$
as $n\to\infty$.  This is the same accuracy as would be expected for sums of
independent indicator random variables, and is thus a remarkably precise result.
However, their proof requires long and delicate analysis of the corresponding
generating functions.  The purpose of this paper is to derive their result
by purely probabilistic arguments, to complement their
result with a distributional approximation in total variation, and to
investigate the quantities~$K_{n,r}$ as well.

The approach that we take begins with the
well--known observation
that, if the fixed value~$n$ were replaced by a Poisson distributed
random number with mean~$n$, then the numbers of balls in the
boxes would be independent Poisson random variables.  Approximations of
the kind to be discussed would then be immediate, from the theory of
sums of independent Bernoulli random variables.  The essence of the
problem lies in the dependence introduced by fixing~$n$.  One way of
relaxing this dependence is to disregard the first few boxes, for which
the result is essentially known, and to use the fact that the number
of balls falling in the remaining boxes is now random.
Indeed, defining $j_n \ge 1$ in such a way that
\eq\label{jn-def}
   p_{j_n-1} \Ge 4n^{-1}\log n\ > \ p_{j_n},
\en
it is immediate that
\[
    \pr[N_j \ge 1\ \mbox{for all}\ j \le j_n - 1] \Ge 1 - \frac n{4\log n}
       \Bigl(1 - \frac {4\log n}n \Bigr)^n \Ge 1 - n^{-3},
\]
so that, except on a set of probability at most~$n^{-3}$, we have
\eq\label{first-ones}
    \sum_{j=1}^{j_n-1} I_j \Eq j_n-1,
\en
where $I_j := I[N_j \ge 1]$.  Furthermore, a simple Poisson approximation argument,
due to Le~Cam~(1960) and Michel~(1988), can now be used to get a sharp description
of the distribution of the remaining elements in the sum $K_n := \sji I_j$, since
\[
   \dtv(\law(N_j,\,j\ge j_n),\law(L_j,\,j\ge j_n)) \Le P_n \ :=\  \sjjn p_j,
\]
where $(L_j,\,j\ge j_n)$ are {\it independent\/} Poisson random variables
with means $\ex L_j = np_j$: see Barbour \& Gnedin~(2009, Section~2).
This means that the random sequences $(I_j,\,j\ge j_n)$ and $(I[L_j\ge1],\,
j\ge j_n)$ can be constructed to be identical, except on a set of
probability at most~$P_n$, so that, except on a set of probability at
most $n^{-3}+P_n$, the distribution of~$K_n$ agrees with that of a sum
of independent indicators, the first $j_n-1$ of which are equal to~$1$.
Hence a discretized central limit theorem and uniform approximation of
point probabilities follow, using $\nn(\m_n,\s_n^2)$ as basis, with accuracies
$O(\s_n^{-1} + n^{-3}+P_n)$ and $O(\s_n^{-2} + n^{-3}+P_n)$ respectively, and
analogous results are also true for the statistics~$K_{n,r}$.

The drawback to this very simple approach is that it need not be the case
that, for instance, $P_n = O(\s_n^{-2})$.  For example, Karlin's case of
regular variation allows the possibility of having $\s_n^2 \asymp n^\b$,
for any given~$\b$, $0 < \b < 1$.  In such cases, $P_n \asymp (n^{-1}\log n)^{1-\b}$,
so that $P_n = O(\s_n^{-2})$ is {\it not\/} true if $\b>1/2$, and
$P_n = O(\s_n^{-1})$ is not true if $\b>2/3$.  To get the result of
Hwang \& Janson~(2008),
we in general need something sharper.

Our approach involves a technique analogous to that above, discarding a set of indices for
which the outcome is essentially known, and using the randomness in the remainder.
Foregoing the total independence of the above scheme, which costs too much to
achieve, we instead construct a {\it conditionally independent\/} sequence
of Binomial random variables within the problem, and use these to provide
the necessary refinement.  The way in which this can be done is described
in R\"ollin~(2005).  There, and in this paper too, we use translations of Poisson
distributions as approximations, instead of discretized normal distributions, though, to the
accuracies being considered, they are equivalent:  the translated Poisson
distribution $\TP(\m,\s^2)$ is defined to be that of the sum of an
{\it integer\/}~$a$ and a Poisson
$\Po(\l)$--distributed random variable, with $\l$ and~$a$ so chosen that
$a+\l = \m$ and $\s^2 \le \l < \s^2 + 1$.

Using this approach, we are able to prove the following two theorems.
We use $\dtv$ to denote the total variation distance between distributions:
\[
    \dtv(P,Q)\ :=\ \sup_A|P(A)-Q(A)|,
\]
and $\dloc$ to denote the local distance
(point metric) between distributions on the integers:
\[
   \dloc(P,Q)\ :=\ \sup_{j\in\Z}|P\{j\} - Q\{j\}|.
\]
We define $j_0$ so that
\[
    \sum_{j\ge j_0-1} p_j \Ge 1/2 \ >\  \sum_{j\ge j_0} p_j =: P_0,
\]
and let $n_0 \ge 3$ be such that $j_n$, defined in~\Ref{jn-def}, satisfies $j_n \ge j_0$
for all $n\ge n_0$, and also that $n_0/\log^2 n_0 \ge 16/P_0$.

\begin{theorem}\label{Kn-Th}
If $\m_n := \ex K_n$ and $\s^2_n := \var K_n$, then
\eqs
   \dtv(\law(K_n),\TP(\m_n,\s^2_n)) &=& O(\s_n^{-1});\\
   \dloc(\law(K_n),\TP(\m_n,\s^2_n)) &=& O(\s_n^{-2}),
\ens
uniformly in $n\ge n_0$.
\end{theorem}

\begin{theorem}\label{Knr-Th}
For $r\ge1$, setting $\m_{n,r} := \ex K_{n,r}$ and $\s^2_{n,r} := \var K_{n,r}$,
we have
\eqs
   \dtv(\law(K_{n,r}),\TP(\m_{n,r},\s^2_{n,r})) &=& O(\s_{n,r}^{-1});\\
   \dloc(\law(K_{n,r}),\TP(\m_{n,r},\s^2_{n,r})) &=& O(\s_{n,r}^{-2}),
\ens
uniformly in $n\ge \max\{n_0 ,e^{r/4},2r\}$.
\end{theorem}

R\"ollin's theorem and our construction are set out in Section~2,
together with the general scheme of the proofs.  The details for the
two theorems are then given in Sections 3 and~4.  Some useful
technical results are collected in the appendix.

\section{The basic method}
\setcounter{equation}{0}
We begin with the following theorem from R\"ollin~(2005).  Let~$W$ be an
integer valued random variable, with mean~$\m$ and variance~$\s^2$, and
let~$M$ be some random element.  Define
\eqa
   \m_M  &:=& \ex(W\giv M);\qquad \s^2_M \ :=\ \var(W\giv M); \qquad \t^2\ :=\ \var(\m_M);\non\\
   \r^2  &:=& \ex(\s^2_M);\qquad \n^2\ :=\ \var(\s^2_M);\qquad U\ :=\ \t^{-1}(\m_M - \m).
   \label{Ro-defs}
\ena
Of course, $\s^2 = \r^2 + \t^2$.

\begin{theorem}\label{Roellin}
Suppose that, for some $\e > 0$,
\eq\label{Stein-eq}
     |\ex\{f'(U) - Uf(U)\}| \Le \e \fdd
\en
for all bounded functions~$f$ with bounded second derivative.  Then there
exist universal constants $R_1$ and~$R_2$ such that
\eqa
   &&\dtv(\law(W),\TP(\m,\s^2)) \non\\
   &&\qquad\le\ \ex\{\dtv(\law(W\giv M),\TP(\m_M,\s^2_M))\}
      + R_1  \frac1\r\Blb 1 + \frac\n{\r} + \frac{\e\t^3}{\s^2}\Brb;\non\\
   &&\dloc(\law(W),\TP(\m,\s^2)) \non\\
   &&\qquad\le\ \ex\{\dloc(\law(W\giv M),\TP(\m_M,\s^2_M))\}
      + R_2 \frac1{\r^2}\Blb 1 + \frac{\n^2}{\r^2} + \frac{\e\t^3}{\s^2}\Brb.\non
\ena
\end{theorem}

\nin Values of the constants are given in R\"ollin~(2005).
Note that~\Ref{Stein-eq}
is exactly what has to be established for the simplest smooth metric standard normal
approximation to~$\law(U)$, using Stein's method.  For~$U$ a sum of independent random
variables, $\e$ would typically be the Lyapounov ratio, and thus the quantity
$\s^{-2}\t^3\e$ would be bounded by an average of the ratios of third to second
moments of the summands.

The theorem is useful provided that $\law(W\giv M)$ is such that it is well
approximated for each value of~$M$ by the translated Poisson distribution
with its mean and variance as parameters.  This is the case, for instance,
for sums of independent Bernoulli random variables, as well as for many sums
of independent integer valued random variables, as noted in R\"ollin~(2005).
Here is the result that we shall use in what follows.

\begin{theorem}\label{Bernoulli}
Suppose  that $\law(W\giv M)$ is the distribution of a sum $\sji
I_j(M)$ of {\it independent\/} Bernoulli random variables with
probabilities $p_j(M)$ such that $\m_M := \sji p_j(M) < \infty$
a.s.; write $\s^2_M := \sji p_j(M)(1-p_j(M))$, $\r^2 := \ex(\s^2_M)$
and $\n^2 := \var(\s^2_M)$. Suppose that $\n^2  \le C\r^2$ for some
$C < \infty$.  Then there exists universal constants $C_1$ and~$C_2$
such that
\eqs
   \ex\{\dtv(\law(W\giv M),\TP(\m_M,\s^2_M))\}
        &\le&  \frac{4C}{\r^2} + \frac{C_1\sqrt2}{\r};\\
   \ex\{\dloc(\law(W\giv M),\TP(\m_M,\s^2_M))\}
        &\le&  \frac{4C+2C_2}{\r^2}.
\ens
\end{theorem}

\proof
Bounds of the form
\eqa
   \dtv(\law(W\giv M),\TP(\m_M,\s^2_M)) &\le&  \min\{C_1 \s^{-1}_M, 1\};\non\\
   \dloc(\law(W\giv M),\TP(\m_M,\s^2_M)) &\le&  \min\{C_2 \s^{-2}_M, 1\},\label{Be-1}
\ena
are given in Barbour~(2009; Theorems 6.2 and~6.3), with $C_1 = 4$ and $C_2 = 280$.
The former follows as in Barbour \& \v Cekanavi\v cius~(2002, Theorem 3.1), and
similar techniques can be used to establish the latter; see also R\"ollin~(2005).
Then, by Chebyshev's inequality, $\pr[\s^2_M < \half\r^2] \le 4C/\r^2$.
The bounds follow by taking expectations in~\Ref{Be-1}.
\ep

\bsk
We now need to find a suitable collection of conditionally independent
Bernoulli random variables.  To do so, we start by observing, as before,
that it is enough to consider indices $j\ge j_n$ in the sums, so we need
only consider the distribution of $(N_j,\,j\ge j_n)$.  We realize these
random variables in two stages: first, we realize $M := (M_j,\, j\ge j_0)$ by
throwing~$n$ balls independently into the boxes with indices $j\ge j_0$,
with probability $p_j/P_0$ for box~$j$, and then `thinning' them
independently with retention probability~$P_0$, so that, conditionally on~$M$,
the~$(N_j,\,j\ge j_0)$ are independent, with $N_j \sim \Bi(M_j,P_0)$.
With this construction, it remains to evaluate the quantities appearing
in R\"ollin's theorem, and to check that we have the right result.
More specifically, we need to check that, for some constants $C,C', C''$,
\eq\label{conditions}
   \mbox{(i) \,} \n^2 \Le C\r^2; \quad \mbox{(ii) \,}\r^2 \Ge C'\s^2, \quad\mbox{and}\quad
     \mbox{(iii) \,} \e \Le C'' \t^{-3}\s^2,
\en
uniformly in the stated ranges of~$n$, for the random variables $W_n :=
\sjjn I[N_j \ge 1]$ and $W_{n,r} := \sjjn I[N_j=r]$, $r\ge1$.  Theorems
\ref{Kn-Th} and~\ref{Knr-Th} will then follow directly from Theorems \ref{Roellin}
and~\ref{Bernoulli}.

The first two inequalities in~\Ref{conditions} cause no great problems,
since they involve only variance calculations, though care has to be
taken with the correlations in Theorem~\ref{Knr-Th}, because the summands
in
\[
    \m_M\ :=\ \sjjn {M_j \choose r} P_0^r (1-P_0)^{M_j-r}
\]
are not monotone functions of the (negatively associated)~$M_j$.  The main effort
is required in evaluating~$\e$ for the third inequality.  We now sketch the
structure of this argument, leaving the details to the next two sections.

Take~$z(l)$, $l\ge0$, to be either
$\Bi(l,P_0)\{[1,\infty)\}$ or $\Bi(l,P_0)\{r\}$, as appropriate, (zero
if $l=0$). Then define the quantity~$U$ that we wish to address by $U := \sjjn Y_j$, where
\eq\label{setting}
   \z_j\ :=\ \ex(z(M_j)), \quad y_j(l)\ :=\ z(l) - \z_j \quad\mbox{and}\quad
     Y_j \ :=\ \t^{-1}y_j(M_j).
\en
Thus~$U$ is a sum of mean zero, weakly dependent random variables.  In order to
approach~\Ref{Stein-eq}, we begin by writing
\eq\label{UfU}
  \ex\{Uf(U)\} \Eq \sjjn \ex\{Y_jf(U)\} \Eq
   \t^{-1}\sjjn\slo q_j(l)y_j(l)\ex\{f(\wjnl + \tyl)\},
\en
where $q_j(l) := \pr[M_j=l]$ and
\eq\label{U-def}
   U_j\um\ :=\ \t^{-1}\ssjnj y_s(M_{js}\um),
\en
and where
\eq\label{multinom}
    M_{j\cdot}\um \ :=\ (M_{js},\,s\ge j_n,s\ne j)\
           \sim\  \MN(m;\,(p_s/P_{0j},\,s\ge j_n,s\ne j))
\en
is distributed as~$m$ balls thrown independently into the boxes with indices
$(s\ge j_n, s\ne j)$ with probabilities $(p_s/P_{0j},\,s\ge j_n,s\ne j)$,
with $P_{0j} := P_0 - p_j \ge 3P_0/4$.  We need to show that the expression
in~\Ref{UfU} is close to $\ex\{f'(U)\}$.

As a first step, we use Taylor development to discard all but the constant and
linear terms in $\ex\{f(\wjnl + \tyl)\}$, establishing that
\eqa
  &(1)&
   \Bigl|\tssqy \{\ex f(\wjnl + \tyl) - \ex f(\wjnl) - \tyl \ex f'(\wjnl)\}\Bigr|\non\\
      &&\qquad   \Le k_1 \stt \fdd. \label{1}
\ena
The next step is to remove the $l$-dependence in the constant term, replacing
$\wjnl$ by $\wjn$.  To make the computations, we realize $\wjnl$ and~$\wjn$
on the same probability space by writing $\mjn = \mjnl + \zjl$, where $\mjnl$
and~$\zjl$ are independent, and distributed as $M_{j\cdot}\um$ in~\Ref{multinom},
with $m = n-l$ and $m=l$, respectively; and then defining $\wjnl$ and~$\wjn$
as before, using~\Ref{U-def}.  Using this representation, we then show that
\eqa
  &(2)&  \Bigl|\tssqy \{\ex f(\wjnl) - \ex f(\wjn)
              - \ex [f'(\wjnl)(\wjnl - \wjn)]\}\Bigr| \non\\
         &&\qquad   \Le k_2 \stt \fdd. \phantom{XXXXXXXXXXXXXXXXXXXXXXXXXXXXX}\label{2}
\ena
Although this has introduced a further term $\ex [f'(\wjnl)(\wjnl - \wjn)]$ involving~$l$,
there is simplification because $\ex f(\wjn)$ is multiplied by $\slo q_j(l)y_j(l) = \ex Y_j
= 0$, and hence drops out.

We now simplify what is left by showing that
\eqa
  &(3)& \Bigl|\tssqy \{\ex [f'(\wjnl)(\wjnl - \wjn)]
              - \ex [f'(\wjn)] \ex(\wjnl - \wjn)\}\Bigr|  \non\\
       &&\qquad     \Le k_3 \stt \fdd. \label{3}
\ena
As a result of this, the quantity $\ex f(\wjnl)$ in~(1) has been replaced by a
multiple of $\ex f'(\wjn)$, with errors of the desired order, which is a useful
step in approaching the intended
goal of $\ex f'(U)$.  There is also the quantity $\ex f'(\wjnl)$ appearing in~(1),
but this is easily reduced to one involving only $\ex f'(\wjn)$, too:
\eq\label{4}
  (4)\quad \Bigl|\tssq y_j^2(l) \{\ex f'(\wjnl) - \ex f'(\wjn) \}\Bigr|
         \Le k_4 \stt \fdd.\phantom{XXXXX}
\en
At this point, we have thus established that
\eq\label{halfway}
  \Bigl| \ex Uf(U) - \t^{-2}\sjjn \k_j \ex f'(\wjn)\Bigr|
     \Le (k_1+k_2+k_3+k_4) \stt\fdd,
\en
with
\eq\label{c-def}
  \k_j \ :=\ \sqy\{y_j(l) - \t \ex(\wjn-\wjnl)\},
\en
and, for example by taking $f(x)=x$,
\[
   1 \Eq \ex U^2 \Eq \t^{-2}\sjjn \k_j.
\]

In parallel with the above reduction starting from~\Ref{UfU}, we now start with
\eqa
  \ex f'(U) &=& \t^{-2}\sjjn \k_j \ex f'(U) % \non\\  &=&
   \Eq \t^{-2}\sjjn \k_j \slo q_j(l) \ex f'(\wjnl + \tyl),\phantom{XXXX}
      \label{fdU}
\ena
and make two rather simpler steps, first proving that
\eq\label{5}
  (5)\ \, \Bigl|\t^{-2}\sjjn \k_j \slo q_j(l) \{\ex f'(\wjnl + \tyl) - \ex f'(\wjnl) \}\Bigr|
         \Le k_5 \stt \fdd,
\en
and then that
\eq\label{6}
  (6)\quad  \Bigl|\t^{-2}\sjjn \k_j \slo q_j(l) \{\ex f'(\wjnl) - \ex f'(\wjn) \}\Bigr|
         \Le k_6 \stt \fdd. \phantom{XXXXXX}
\en
Putting these two into~\Ref{fdU}, it follows that
\eq\label{fulltime}
  \Bigl| \ex f'(U) - \t^{-2}\sjjn \k_j \ex f'(\wjn) \Bigr| \Le (k_5+k_6)\stt\fdd,
\en
and combining this with~\Ref{halfway} yields
\eq\label{job-done}
  |\ex\{f'(U) - Uf(U)\}| \Le \e\fdd,
\en
with $\s^{-2}\t^3\e \le \sum_{t=1}^6 k_t$ bounded, as required.

\section{The argument for~$K_n$}
\setcounter{equation}{0}
We begin by noting, for future reference, that we have
\eqa
    \bp_n &:=& \max_{j\ge j_n} p_j \Le 4n^{-1}\log n \Le P_0/4 \Le 1/8; \non\\
    n\bp_n^2 &\le& 16 n^{-1}\log^2n \Le P_0, \label{bp-bnd}
\ena
whenever $n\ge n_0$, and that $\b := (1-P_0/2) \ge 3/4$.
We use $c$ and~$c'$ to denote generic universal constants, not
depending on~$n$ or the $p_j$'s.

For~$K_n$, we have $\law(W_n \giv M)$ that of a sum of indicator random
variables $I_j(M)$, $j\ge j_n$, with probabilities
\[
    \{1 - (1-P_0)^{M_j}\}\ =:\ z(M_j);
\]
recall~\Ref{setting}.  Hence $\s^2_M = \sjjn z(M_j)(1-z(M_j))$, and
\[
   \r^2 \Eq \ex \s^2_M \Eq \sjjn \ex\{(1-P_0)^{M_j} - (1-P_0)^{2M_j}\}.
\]
Applying Lemma~\ref{idmom}\,(iv) with $x = \sqrt{1-P_0}$,
and using the fact that $n\bp_n^2 \le P_0$, now immediately gives the
lower bound
\eqa
   \r^2 &\ge& c_\r \sjjn e^{-np_j}\min\{1,np_j\}, \label{rho2}
\ena
where $c_\r = c(\sqrt{1-P_0})e^{-2P_0}$, and $c(\cdot)$ is as in Lemma~\ref{idmom}.
On the other hand, because the~$N_j$ are negatively associated,
\[
   \s^2 \Le \sjjn \var I[N_j \ge 1] \Eq \sjjn \{1-(1-p_j)^n\}(1-p_j)^n
    \Le \sjjn e^{-np_j}\min\{1,np_j\}.
\]
It thus follows that $\r^2 \ge c_\r\s^2$, % with $C = \min\{(1 - e^{-c}),ce^{-c}\}$,
establishing~\Ref{conditions}\,(ii).

For $\n^2 = \var \s^2_M$, we note that $\s^2_M$ is the difference of the random
variables $s_1(M) := \sjjn (1-P_0)^{M_j}$ and $s_2(M) := \sjjn (1-P_0)^{2M_j}$, so that
$\n^2 \le 2(\var s_1(M) + \var s_2(M))$.  Since $(1-P_0)^l$ is decreasing
in~$l$, we can use the negative association of the $M_j$'s to upper bound
the variances:
\[
   \var s_1(M) \Le \sjjn \var \{(1-P_0)^{M_j}\};
             \qquad \var s_2(M) \Le \sjjn \var\{(1-P_0)^{2M_j}\}.
\]
Now both of these quantities can be bounded by using Lemma~\ref{idmom}\,(iv):
\[
   \var \{(1-P_0)^{M_j}\} \Le e^{-2\b np_j}\min\{1,2\b np_j\},
\]
and
\[
   \var \{(1-P_0)^{2M_j}\} \Le e^{-2\b' np_j}\min\{1,2\b' np_j\},
\]
with $\b' := 4 - 6P_0 + 4P_0^2 - P_0^3$.  Thus $\r^{-2}\n^2$ is
uniformly bounded, establishing~\Ref{conditions}\,(i).  It thus remains to prove that
$\e \le C''\t^{-3}\s^2$ for some constant~$C''$, and we are finished.
To do this, we successively verify the inequalities (1) -- (6) of
Section~2.

To establish inequality~(1), we note that its left hand side is
bounded by
\eq\label{1-zero}
    \half\t^{-3} \sjjn\slo q_j(n)|y_j(l)|^3 \fdd.
\en
Now $|y_j(l)| \le 1$, and
\[
    \sq y_j^2(l) \Eq \ex\{(1-P_0)^{2M_j}\} - \{\ex(1-P_0)^{M_j}\}^2,
\]
with $M_j \sim \Bi(n,p_j/P_0)$. From Lemma~\ref{idmom}\,(iv)
with $x=1-P_0$, it follows that
\eqa
   \sq y_j^2(l) &\le& e^{-2\b np_j}\min\{1,2\b np_j\}. \label{2.35}
\ena
Hence, from Lemma~\ref{sum-bds}\,(i),
\[
   \t^{-3}\sq |y_j(l)|^3 \Le \t^{-3}\sjjn np_j e^{-2\b np_j} \Le K_0^{(2\b-1)}\stt.
\]
By~\Ref{1-zero}, this proves~(1) with $k_1 = K_0^{(2\b-1)}$.

For inequality~(2), we have
\eq\label{2.41}
   |\ex\{f(\wjn) - f(\wjnl) - f'(\wjnl)(\wjn-\wjnl)\}|
      \Le \half\fdd \ex\{(\wjn-\wjnl)^2\}.
\en
Now
\[
  \t^{2}\ex\{(\wjn-\wjnl)^2\} \Le \ex\Blb \Bl \ssjnj \zjsl P_0(1-P_0)^{\mjsnl}\Br^2 \Brb,
\]
and the collections of random variables $(\zjsl,\,s\ge j_n)$ and $((1-P_0)^{\mjsnl},\,s\ge j_n)$
are independent, and each is composed of negatively correlated elements.  Hence
\eqs
  \lefteqn{\t^{2}\ex\{(\wjn-\wjnl)^2\}}\\
  && \Le P_0^2\Bl \ssjnj \ex \zjsl\, \ex\Bigl\{(1-P_0)^{\mjsnl}\Bigr\}\Br^2
    + P_0^2 \ssjnj \ex\{(\zjsl)^2\}\ex\Bigl\{(1-P_0)^{2\mjsnl}\Bigr\}.
\ens
Now routine calculation gives
\eqs
  &&P_0\,\ex\zjsl \Le lP_0p_s/P_{0j} \Le 2lp_s; \qquad P_0^2\,\ex\{(\zjsl)^2\} \Le 2lp_s(1+2lp_s);\\
  &&\ex\Bigl\{(1-P_0)^{\mjsnl}\Bigr\} \Le e^{-(n-l)p_s};\qquad
            \ex\Bigl\{(1-P_0)^{2\mjsnl}\Bigr\} \Le e^{-2\b(n-l)p_s},
\ens
and hence, with crude simplifications,
\eq\label{2.42}
  \t^{2}\ex\{(\wjn-\wjnl)^2\} \Le 10l^2e^{l\d_n}\ssjn p_s e^{-2\b np_s}
    \Le c l^2 e^{l\d_n}  n^{-1}\s^2 ,
\en
this last using~\Ref{rho2} and Lemma~\ref{sum-bds}\,(i), where $\d_n := 2\bp_n$
and $c=10(K(2\b-1)/c_\r)$.
Hence, putting \Ref{2.41} and~\Ref{2.42} into~(2), we obtain the bound
\eqs
   \lefteqn{\frac{c}2\fdd \t^{-3}\ssq |y_j(l)|l^2 e^{l\d_n}n^{-1}\s^2}\\
     &&\Le c'\t^{-3} \s^2\fdd \exp\{\d_n(3 + n\bp_n e/P_0)\} \sjjn e^{-np_j} p_j(1 + np_j),
\ens
by Lemma~\ref{idmom}\,(ii) and~(iii), and this is uniformly of order~$\t^{-3}\s^2\fdd$ in
the stated range of~$n$, because
\[
   \sjjn p_j(1+np_j)e^{-np_j} \Le P_n(1+e^{-1}) \quad\mbox{and}\quad
     \d_n + n\d_n\bp_n \Le 5P_0/4.
\]
This establishes inequality~(2).

For inequality~(3), we begin by writing
\eqa
  \lefteqn{\ex\{(\wjnl - \wjn)f'(\wjnl)\} }\non\\
  &&\Eq \ex\{[\ex(\wjnl - \wjn \giv \mjnl) - \ex(\wjnl - \wjn)](f'(\wjnl)-f'(\ex\wjnl))\} \non\\
   &&\mbox{}\qquad\qquad   - \ex(\wjn - \wjnl)\ex f'(\wjnl);  \label{2.50}
\ena
note that introducing $f'(\ex\wjnl)$ changes nothing, since it is multiplied by
a quantity with mean zero.  The first term we bound by
\eq\label{var-bnds}
   \fdd\,\sqrt{\var[\ex(\wjnl - \wjn \giv \mjnl)]}\,\sqrt{\var\wjnl}.
\en
Since
\eq\label{2.52}
   \t \ex(\wjnl - \wjn \giv \mjnl) \Eq \ssjnj (1-P_0)^{\mjsnl} \{1 - (1-p_sP_0/P_{0j})^l\},
\en
and since the $(\mjsnl,\, s\ge j_n)$ are negatively associated, it follows that
\eqs
   \t^2 \var[\ex(\wjnl - \wjn \giv \mjnl)] &\le& 4l^2\ssjnj p_s^2e^{-2\b(n-l)p_s}\\
     &\le& 4l^2 e^{l\d_n}n^{-1}/(2\b e) \Eq cl^2e^{l\d_n}n^{-1},
\ens
for a suitable~$c$.
In much the same way, and using Lemma~\ref{idmom}\,(iv), we have
\[
  \t^2\var\wjnl \Le \ssjnj \var\{(1-P_0)^{\mjsnl}\}
     \Le 2\frac{P_0}{P_{0j}}\ssjnj np_s e^{-2\b(n-l)p_s} \Le ce^{l\d_n} \s^2.
\]
Hence the first term in~\Ref{2.50} is bounded by
\eq\label{2.55}
  c\t^{-2}\fdd\,le^{l\d_n}n^{-1/2}\s,
\en
for a suitable~$c$.  For the second, we replace $\ex f'(\wjnl)$ by $\ex f'(\wjn)$:
\eq\label{2.60}
  |\ex(\wjn - \wjnl)\{\ex f'(\wjnl) - \ex f'(\wjn)\}| \Le \fdd \ex\{(\wjn-\wjnl)^2\},
\en
which is at most $c\t^{-2}\fdd l^2 e^{l\d_n} n^{-1}\s^2$.  Putting these bounds
into~\Ref{2.50}, it follows that the left hand side in~(3) is at most
\eqa
  \lefteqn{
     c\t^{-3}\fdd\ssq |y_j(l)| e^{l\d_n}\{ln^{-1/2}\s + l^2n^{-1}\s^2\}
  } \non\\ && \Le
%  &\le&
    c'\t^{-3}\fdd\Blb n^{-1/2}\s \sjjn np_j e^{-np_j} + \s^2\Brb, \label{2.63}
\ena
by using Lemma~\ref{idmom}\,(ii) and~(iii), for suitable constants $c$ and~$c'$.  But now
\[
  \sjjn np_j e^{-np_j} \Le \sqrt{K'n\s^2},
\]
by Lemma~\ref{sum-bds}\,(iv), and this, together with~\Ref{2.63}, shows
that~(3) is satisfied. % with a suitable choice of~$k_3$.

For~(4), we use the simple bound
\eq\label{simple-l-bd}
  |\ex f'(\wjnl) - \ex f'(\wjn)| \Le \fdd\,\ex|\wjn-\wjnl| \Le \t^{-1}l\fdd.
\en
This gives a bound for the left hand side of~(4) of
\eqs
  \t^{-3}\fdd \ssq y_j^2(l) l \Le \t^{-3}\fdd \sjjn np_j\{e^{-2np_j} + e^{-2\b np_j}\}
    \Le k_4 \t^{-3}\fdd \s^2,
\ens
by Lemma~\ref{sum-bds}\,(i); and hence we have proved~\Ref{halfway}.

For the remaining two inequalities, we observe that, from~\Ref{c-def}
and~\Ref{2.35},
\eq\label{cj-1}
   \k_j^+ \ :=\ \max\{\k_j,0\} \Le 2\b np_j e^{-2\b np_j},
\en
whereas, from~\Ref{2.52},
\eq\label{cj-2}
  \k_j^- \:=\ |\min\{0,\k_j\}| \Le \sq|y_j(l)|\ssjn 2lp_s e^{-(n-l)p_s}
     \Le cnp_je^{-np_j}\ssjn p_se^{-np_s},
\en
from Lemma~\ref{idmom}\,(ii) and~(iii).
Hence, for inequality~(5), we obtain the bound
\eqa
   \t^{-3}\fdd \sjjn |\k_j| \slq|y_j(l)| &\le& 2 \t^{-3}\fdd \sjjn |\k_j| e^{-np_j} \non\\
  &\le& c \t^{-3}\fdd \sjjn np_j e^{-2np_j} \Le k_5\t^{-3}\s^2\fdd, \phantom{XX}\label{2.72}
\ena
by Lemma~\ref{sum-bds}\,(i), for a suitable~$k_5$.  For inequality~(6), we start from the bound
\eqs
   \lefteqn{ \t^{-2}\fdd \sjjn |\k_j| \slq \ex|\wjn - \wjnl| }\\
    &&\le\ \t^{-3}\fdd \sjjn |\k_j| \slq 2le^{l\d_n}\ssjnj p_s e^{-np_s}  %\\ &&
    \Le c\t^{-3}\fdd \sjjn |\k_j| np_j\ssjn p_s e^{-np_s},
\ens
again from~\Ref{2.52} and Lemma~\ref{idmom}\,(ii), and substituting from
\Ref{cj-1} and~\Ref{cj-2} for~$|\k_j|$ gives at most
\eqa
%  \lefteqn{
    c\t^{-3}\fdd \sjjn (np_j)^2
     \Blb  P_ne^{-2\b np_j} +  e^{-np_j}\Bl\ssjn p_se^{-np_s}\Br^2 \Brb %} \non\\ &&
   \Le k_6 \t^{-3}\fdd\,\s^2, %\phantom{XXXXXXXXXXXXXXXXXX}
   \label{2.73}
\ena
% for suitable choice of~$k_6$,
by Lemma~\ref{sum-bds}\,(i) and~(iv).
Since \Ref{2.72} and~\Ref{2.73} together establish~\Ref{fulltime}, we have completed
the proof of~\Ref{job-done}, and hence of~\Ref{conditions}\,(iii), thus proving Theorem~\ref{Kn-Th}.

\section{The argument for~$K_{n,r}$}
\setcounter{equation}{0}
Fix~$r\ge1$. We now require~$n$ to satisfy $4\log n \ge r-1$ and $n \ge 2r$.
Then, with $p := p_{j_n-1} \ge 4n^{-1}\log n$, we have
\eqs
  &&\sum_{j < j_n} \pr[N_j=r] \Le (j_n-1){n\choose r}p^r(1-p)^{n-r}
      \Le n^r p^{r-1} e^{-(n-r)p} /r!\\
    &&\qquad\Le n^{-3}(4\log n)^{r-1} e^r/r! \Le c(\log n)^{r-1}n^{-3},
\ens
since $x^s e^{-x}$ is decreasing in $x\ge s$ and $4\log n \ge r-1$.
Hence $\sum_{j < j_n} I[N_j=r] = 0$ except on a set of probability
of order~$O(n^{-3}(\log n)^{r-1})$, and we can restrict attention to
$W_{n,r} := \sjjn I[N_j=r]$.  We recall that  $\b := (1-P_0/2) \ge 3/4$, and that
\[
    \bp_n \Le P_0/4 \Le 1/8 \quad\mbox{and}\quad n\bp_n^2  \Le P_0,
\]
whenever $n\ge n_0$.
% We now also require~$n$ to satisfy $\min\{4\log n,n/2\} \ge r$.
The generic constants $c$ and~$c'$ are now allowed to depend on~$r$.

For~$K_{n,r}$, the distribution $\law(W_{n,r} \giv M)$ is that of a sum of indicator random
variables $I_j(M)$, $j\ge j_n$, with probabilities
\[
    {M_j \choose r}P_0^r(1-P_0)^{M_j}\ =:\ z(M_j);
\]
recall~\Ref{setting}.  The argument now runs much as before, but
is complicated by the fact that $z(\cdot)$ is not monotonic in~$l$.
First, we have $\m = \sjjn \ex z(M_j) = \sjjn \z_j$, with
$\z_j := \Bi(n,p_j)\{r\}$, whence, defining
\ignore{
\[
   \z_j\ :=\ {n\choose r}p_j^r(1-p_j)^{n-r},
\]
and it follows easily that, setting
}
\[
    \hm_r \ :=\ \sjjn\frac{(np_j)^r e^{-np_j}}{r!},
\]
it easily follows that
\eq\label{mu-bd}
  \exp\{-n\bp_n^2 - n^{-1}r^2\}
     \Le \m /\hm_r
        \Le e^{r\bp_n},
\en
for $n\ge 2r$,
with both lower and upper estimates uniformly bounded away from zero
and infinity in the chosen range of~$n$: hence $\m$ and~$\hm_r$ are
uniformly of the same order.

Now
\eq\label{s2M-def-r}
    \s^2_M \Eq  \sjjn z(M_j)(1-z(M_j)) \Ge \sjjn z(M_j)(1-z_r),
\en
where $z_r := \max_{l\ge r} {l\choose r} P_0^r(1-P_0)^{l-r} < 1$, and
hence
\eq\label{rho-lower}
   \r^2 \Eq \ex \s^2_M \Ge \m (1-z_r).
\en
For
\[
   \s^2 \Eq \var W_n \Eq \sjjn \ssjn \{\pr[N_j=N_s=r] - \pr[N_j=r]\pr[N_s=r]\},
\]
we use Lemma~\ref{cov-bd} to give
\[
   \pr[N_j=N_s=r] - \pr[N_j=r]\pr[N_s=r] \Le 2er(p_j+p_s)e^{4r\bp_n}\pr[N_j=r]\pr[N_s=r],
   \quad j\ne s,
\]
and adding over $j$ and~$s$ gives an upper bound of at most
\[
  c\sjjn p_j(np_j)^r e^{-np_j} \ssjn (np_s)^r e^{-np_s}
  \Le c'P_n\hm_r.
\]
For $j=s$, the total contribution to the variance is at most $\sjjn \pr[N_j=r] = \m$.
Hence, and from~\Ref{rho-lower}, we have
\eq\label{asymp-similar}
   \s^2\ \asymp\ \r^2\ \asymp\ \m\ \asymp\ \hm_r,
\en
where the implied constants are universal for each~$r$.
This shows also that~\Ref{conditions}\,(ii) holds.

For~\Ref{conditions}\,(i), we take
\[
   \n^2 \ :=\ \var(\s^2_M) \Eq \var\Bl \sjjn z(M_j)(1-z(M_j)) \Br,
\]
to which we can apply Lemma~\ref{cov-bd}, noting that
$0 \le z(l)(1-z(l)) \le {l\choose r}P_0^r(1-P_0)^{l-r}$.
For $j\ne s$, this gives
\[
  \cov\{z(M_j)(1-z(M_j)),z(M_s)(1-z(M_s))\}
     \Le c(p_j+p_s)(n(p_j+p_s)+2r)(np_j)^r(np_s)^r e^{-n(p_j+p_s)},
\]
by Lemma~\ref{mdmom}. Adding over $j$ and~$s$, this gives at most
\eq\label{3.9}
    c'\Blb\sjjn p_j(np_j + 2r)(np_j)^r e^{-np_j} \ssjn (np_s)^r e^{-np_s}
         + \sjjn p_j(np_j)^r e^{-np_j} \ssjn (np_s)^{r+1} e^{-np_s} \Brb,
\en
and this is at most $cP_n\hm_r + K_{11}P_n\hm_r$, by Lemma~\ref{sum-bds}\,(iii)
and~(iv).  The terms with $j=s$ give at most
\eqa
   \sjjn \ex\{z^2(M_j)\} &\le& \frac{P_0^{2r}}{(r!)^2}
          \ex\Bigl\{ [(M_j)_{(2r)} + (2r)\dr(M_j)\dr](1-P_0)^{2(M_j-r)} \Bigr\} \non\\
    &\le& c\{(np_j)^{2r} + (np_j)^r\} e^{-2\b(n-r)p_j},  \label{3.10}
\ena
by Lemma~\ref{idmom}, and because $l^2\dr \le {2r\choose r}l_{(2r)} + (2r)\dr l\dr$.
Adding over~$j$, this gives at most a contribution of $c\hm_r$, by Lemma~\ref{sum-bds}.
Thus we have shown that $\n^2 \le c\s^2$, and~\Ref{conditions}\,(i)
is satisfied.  It thus remains to show that $\e \le c\t^{-3}\s^2$, and the proof
is accomplished.

To establish inequality~(1), we once again observe that $|y_j(l)|
:= |z(l) - \ex z(M_j)| \le 1$, and hence, recalling~\Ref{1-zero}, that
\[
    \half\t^{-3}\fdd\, \sjjn\ex|y_j(M_j)|^3
      \Le \t^{-3}\fdd\,\sjjn \ex z^2(M_j) \Le c\t^{-3}\fdd\,\hm_r,
\]
as for~\Ref{3.10}; so~(1) holds, as required.

For~(2), we recall~\Ref{2.41}. We then note that, for $u\ge r$,
\eq\label{3.15}
  |z(u+t) - z(u)| \Eq P_0^r\Bigl|{u\choose r}(1-P_0)^{u-r}
            - {u+t \choose r}(1-P_0)^{u+t-r}\Bigr| \Le c{u\choose r}(1-P_0)^u,
\en
for~$c$ a universal constant.  From this, it follows that
\eqa
  \lefteqn{\t|\wjn - \wjnl|} \label{3.16}\\
  && \Le \ssjnj\Blb cI[\zjsl \ge 1]{\mjsnl \choose r}(1-P_0)^{\mjsnl}
     + \sum_{u=0}^{r-1} I[\zjsl \ge r-u]\,I[\mjsnl=u]\Brb. \non
\ena
Since $(x_1+\cdots+x_r)^2 \le r(x_1^2 + \cdots + x_r^2)$, we can bound
$\t^2\ex(\wjn - \wjnl)^2$ by considering the~$r$ different sums separately.

First, for
\[
    \ex \Blb\Bl \ssjnj I[\zjsl \ge 1]{\mjsnl \choose r}(1-P_0)^{\mjsnl}\Br^2 \Brb,
\]
using the independence of $\zjl$ and~$\mjnl$ and Lemma~\ref{mdmom},
and with $\d_n = 2\bp_n$ as before, the off-diagonal terms give at most
\eqs
   c\ssjn \stjn (l^2 p_sp_t)(np_s)^r(np_t)^r e^{-n(p_s+p_t)} e^{2\d_n(2r+l)}
   \Le c'l^2 e^{2l\d_n}\,n^{-1}P_n\hm_r,
\ens
the last line using Lemma~\ref{sum-bds}\,(v).  The terms with $j=s$
then contribute at most
\[
   c\ssjn lp_s (np_s)^r \{1 + (np_s)^r\}e^{-2\b np_s} e^{2l\d_n}
      \Le c'l e^{2l\d_n} n^{-1}\hm_r,
\]
using Lemma~\ref{sum-bds}\,(ii).  The contribution to $\t^2\ex(\wjn - \wjnl)^2$
from this first sum is thus no more than $cl^2 e^{2l\d_n} n^{-1}\hm_r$

For $0\le u \le r-1$, we need to find similar bounds for
\[
  \ex \Blb\Bl \ssjnj I[\zjsl \ge r-u]I[\mjsnl =u]\Br^2 \Brb.
\]
Here, the off-diagonal terms contribute at most
\eqs
  c\ssjn \stjn (l^{2(r-u)} (p_sp_t)^{r-u}(np_s)^u(np_t)^u
               e^{-n(p_s+p_t)} e^{2\d_n(2u+l)}
   \Le c'(l/n)^{2(r-u)} e^{2l\d_n} n\hm_r, %\phantom{XXXXXXXXXXXXXXX}
\ens
by Lemma~\ref{sum-bds}\,(v), and the diagonal terms give at most
\[
   c\ssjn (lp_s)^{r-u} (np_s)^u e^{-np_s}e^{2\d_n(2u+l)}
      \Le c'(l/n)^{r-u} e^{2l\d_n} \hm_r.
\]
Since, in the above, $u \le r-1$ and $l\le n$, it follows that
\eq\label{w-squared}
    \t^2\ex(\wjn - \wjnl)^2 \Le cl^2 e^{2l\d_n} n^{-1}\hm_r.
\en

Returning to~(2), and once again recalling~\Ref{2.41}, we thus have a bound of
\eqs
  \lefteqn{ \half\fdd\,  \tssq|y_j(l)| \ex(\wjn - \wjnl)^2 % }\\
%    } \\   &&
  \Le c\t^{-3}\fdd\, \frac{\hm_r}n \sjjn \ex\{|y_j(M_j)| M_j^2 e^{2M_j\d_n}\} }\\
  &&\Le c'\t^{-3}\fdd\, \frac{\hm_r}n  \sjjn (np_j)^r(1 + (np_j)^2) e^{-np_j}
%    \\     &&
  \Le c'\hm_r \t^{-3}\fdd\, (K_{r-1} + K_{r+1})P_n,
\ens
from Lemma~\ref{sum-bds}\,(iii), and this completes the proof of~(2).

For inequality~(3), recalling~\Ref{2.50} and \Ref{var-bnds}, we
first need to bound the variance $\var\{\ex(\wjn - \wjnl \giv \mjnl)\}$.  Now
\[
  \t\ex(\wjn - \wjnl \giv \mjnl) \Eq \ssjnj \ex(z(\mjsn) - z(\mjsnl) \giv \mjnl)
    \ =:\ \ssjnj g_s(\mjsnl),
\]
where, from~\Ref{3.15} and the independence of $\zjl$ and~$\mjnl$,
\eq\label{3.27}
   |g_s(t)| \Le \frac{lp_s}{P_{0j}}\,{t\choose r} (1-P_0)^t P_0^r,
\en
but $g_s$ is not non-negative.
From Lemmas \ref{cov-bd} and~\ref{mdmom}, the off-diagonal terms in
the variance $\var\{\sum_{s\ge j_n,\,s\ne j} g_s(\mjsnl)\}$ contribute at most
\[
   cl^2e^{2l\d_n} \ssjn \stjn p_sp_t (np_s)^r(np_t)^r
     \{(p_s+p_t)(1+np_s+np_t) + n^{-1}(1+np_s)(1+np_t) + np_sp_t\}e^{-n(p_s+p_t)},
\]
and, using Lemma~\ref{sum-bds}, this can be bounded by
$cl^2e^{2l\d_n}n^{-2}P_n\hm_r$.  The diagonal terms in turn yield at most
\[
   \ssjnj \var g_s(\mjsnl) \Le cl^2e^{2l\d_n}
       \ssjn p_s^2 (np_s)^r(1 + (np_s)^r) e^{-2\b np_s}
           \Le c'l^2e^{2l\d_n} n^{-1}P_n,
\]
by Lemma~\ref{sum-bds}\,(iii).  Since also $\hm_r \le cn$, it follows that
\[
   \var\{\ex(\wjn - \wjnl \giv \mjnl)\} \Le c \t^{-2}l^2e^{2l\d_n} n^{-1}P_n.
\]

For $\t^2\var\wjnl$, the considerations are similar but easier, since we
now have
\[
    0 \Le z(t) \le {t\choose r}(1-P_0)^t P_0^r
\]
in place of~\Ref{3.27}, and the contributions from both diagonal and
off-diagonal terms are bounded by $e^{2l\d_n}\hm_r$.  Hence, and
recalling~\Ref{2.50} and~\Ref{var-bnds}, we have arrived at a bound
\eqa
  \lefteqn{|\ex\{[\ex(\wjnl - \wjn \giv \mjnl) - \ex(\wjnl - \wjn)]
                                          (f'(\wjnl)-f'(\ex\wjnl))\}|}\non\\
  &&\Le c\t^{-2}\fdd\,l e^{2l\d_n} \sqrt{\hm_rP_n/n};\phantom{XXXXXXXXXXXXXXXXXXXXXXX}
\ena
the analogue of~\Ref{2.60},
\eq\label{2.60-r}
  |\ex(\wjn - \wjnl)\{\ex f'(\wjnl) - \ex f'(\wjn)\}|
           \Le c\t^{-2}\fdd\,l^2 e^{2l\d_n}n^{-1}\hm_r,
\en
follows directly from~\Ref{w-squared}.  Hence, for~(3), we have
\eqs
  \lefteqn{ \Bigl|\tssqy \{\ex [f'(\wjnl)(\wjnl - \wjn)]
              - \ex [f'(\wjn)] \ex(\wjnl - \wjn)\}\Bigr| }\\
  &&\Le c\t^{-3}\fdd\,\sjjn \ex\{M_j^2|y_j(M_j)| e^{2M_j\d_n}\}(\sqrt{\hm_rP_n/n} + n^{-1}\hm_r)\\
  &&\Le c' \t^{-3}\fdd\,\Blb \sjjn (np_j)^{r+1}(1+np_j)e^{-np_j}\Brb
             (\sqrt{\hm_rP_n/n} + n^{-1}\hm_r),\phantom{XXXXXXX}
\ens
and since
\eq\label{square-bd-r}
   \Blb \sjjn (np_j)^{r+1}(1+np_j)e^{-np_j}\Brb^2 \Le cnP_n\hm_r,
\en
by Lemma~\ref{sum-bds}\,(v), we conclude that inequality~(3) is indeed
satisfied.

For inequality~(4), we use the simple bound in \Ref{simple-l-bd},
obtaining
\eqs
  \lefteqn{\Bigl|\tssq y_j^2(l) \{\ex f'(\wjnl) - \ex f'(\wjn) \}\Bigr|
   \Le \t^{-3}\fdd\, \sjjn\ex \{M_j y_j^2(M_j)\} }\\
   &&\Le  c\t^{-3}\fdd\, \sjjn (np_j)^r(1 + (np_j)^{r+1})e^{-2\b np_j}
   \Le c'\hm_r\t^{-3}\fdd, \phantom{XXXXXXXXX}
\ens
from Lemma~\ref{idmom}\,(iii), in much the same way as for~\Ref{3.10}.
Hence we have now established~\Ref{halfway}.

For (5) and~(6), we need the constants~$\k_j$, for which we now have
the bounds
\[
   \k_j^+ \Le c(np_j)^r(1 + (np_j)^r)e^{-2\b np_j},
\]
from~\Ref{3.10}, and
\eqs
   \k_j^- &\le& c\ex\{M_j|y_j(M_j)|e^{2M_j\d_n}\} \sqrt{\hm_r/n} \\
     &\le& c'(np_j)^r(1+np_j)e^{-np_j} \sqrt{\hm_r/n},
\ens
from~\Ref{w-squared}.
For inequality~(5), this immediately gives a bound of
\[
   c \t^{-3}\fdd \sjjn |\k_j| (np_j)^r e^{-np_j}  \Le c'\hm_r\t^{-3}\fdd,
\]
using Lemma~\ref{sum-bds}\,(ii); for~(6), we obtain the bound
\[
   c \t^{-3}\fdd \sjjn |\k_j| np_j \sqrt{\hm_r/n}  \Le c'\hm_r\t^{-3}\fdd,
\]
where, for the contribution from $\k_j^-$, we again use
Lemma~\ref{sum-bds}\,(v), much as for~\Ref{square-bd-r}.
This completes the proof of~\Ref{fulltime}, and thus of Theorem~\ref{Knr-Th}.

\section{Appendix}
\setcounter{equation}{0}
We collect several useful calculations, the first two of which
need little proof.  We write $m\ds := m(m-1)\ldots(m-s+1)$.

\begin{lemma}\label{idmom}
If $M\sim \Bi(m,p)$, then for any $x>0$ and $0\le s\le m$,
\[
  (i)\quad \ex\{M\ds x^M\} \Eq m\ds (xp)^s(1 + p(x-1))^{m-s}.\phantom{XXXXXX}
\]
In particular, if $x = e^\d$, where $0 \le \d \le \d_0 \le 1$, and if
$(1-P)e^{\d_0} \le 1$, then
\eqs
   (ii)&&  \ex\{M\ds x^M\} \Le (mp)^s \exp\{\d_0(s+mpe)\};\\
   (iii)&&  \ex\{M\ds [(1-P)e^{\d}]^M\} \Le (mp(1-P))^s e^{-(m-s)pP} \exp\{\d_0[s+mpe(1-P)]\}.
\ens
Furthermore, for $0\le x\le 1$ and $p\le 1/2$, we have
\eqs
  &&(iv)\quad   c(x) e^{-2mp^2} \min\{1,mp\} \Le e^{mp(1-x^2)} \{\ex x^{2M} - (\ex x^M)^2\}
    \Le \min\{1,mp(1-x^2)\},
\ens
where $c(x) := \min\{(1-e^{-(1-x)^2}),(1-x)^2 e^{-(1-x)^2})\}$.
\end{lemma}

\proof
We prove only (iv).  From (i), we have
\[
   \ex x^{2M} - (\ex x^M)^2 \Eq \{1-p(1-x^2)\}^m
          \Blb 1 - \Bl 1 - \frac{p(1-p)(1-x)^2}{1-p(1-x^2)} \Br^m\Brb.
\]
The upper bound follows immediately, using the fact that $1-p \le 1-p(1-x^2)$.
The lower bound
\[
    e^{-mp(1-x^2) -2mp^2}\{1 - e^{-mp(1-x)^2}\}
\]
also uses the fact that $p \le 1/2$, and the argument is completed in standard
fashion.
\ep

\bsk
\begin{lemma}\label{mdmom}
Let $(L,M,m-L-M) \sim \MN(m;\,p,q,1-p-q)$ be trinomially distributed. Then
\[
   \ex\{L\du M\dv w^L x^M\} \Eq m_{(u+v)} (wp)^u(xq)^v(1 + p(w-1) + q(x-1))^{m-u-v}.
\]
In particular, if $0 \le w,x \le e^\d$, where $0 \le \d \le \d_0 \le 1$, and if
$(1-P)e^{\d_0} \le 1$, then
\eqs
  &&\ex\{L\du M\dv w^L x^M\} \Le (mp)^u(mq)^v \exp\{\d_0[(u+v) + m(p+q)e]\};\\
  &&\ex\{L\du M\dv [(1-P)e^{\d}]^{L+M}\} \\
     &&\quad\Le (mp(1-P))^u(mq(1-P))^v e^{-(m-u-v)(p+q)P} \exp\{\d_0[(u+v) + m(p+q)e(1-P)]\}.
\ens
\end{lemma}

\bsk
\begin{lemma}\label{cov-bd}
Let $(L,M,m-L-M) \sim \MN(m;\,p,q,1-p-q)$ be trinomial, where $p+q\le\d\le1/4$,
and let the functions $f,g,h,k$ satisfy $0\le f(l) \le h(l)$ and $0 \le g(l) \le k(l)$
for $l\in\Z_+$.  Then
\eqs
  \lefteqn{  \cov(f(L),g(M)) \Le C_1}  \\
   &&:= e(p+q)\{\ex(Lh(L)e^{2L\d}) \ex(k(M)e^{2M\d}) + \ex(h(L)e^{2L\d}) \ex(Mk(M)e^{2M\d})\}.
\ens
If $f$ and~$g$ are not nonnegative, but $|f|$ and~$|g|$ are bounded as above,
then
\eqs
   \cov(f(L),g(M)) &\le& C_1
   + 2m^{-1}\ex(Lh(L)) \ex(Mk(M)) + \frac{4m}3 pq\ex h(L) \ex k(M).
\ens
\end{lemma}

\proof
From the multinomial formulae, we have
\eqa
  \lefteqn{f(u)g(v)\{\pr[L=u,M=v] - \pr[L=u]\pr[M=v]\}}   \non\\
  &&\Eq \frac{f(u)g(v)}{u!v!}p^uq^v
      \{m_{(u+v)}(1-p-q)^{m-u-v} - m\du m\dv (1-p)^{m-u}(1-q)^{m-v}\} \non \\
  &&\Le f(u)g(v)\pr[L=u]\pr[M=v]\{(1-p-q)^{-(u+v)} - 1\} \label{A.10}\\
  &&\Le h(u)k(v)\pr[L=u]\pr[M=v](p+q)(u+v)\exp\{2(p+q)(u+v+1)\}, \non
\ena
where the last inequality uses $p+q \le 1/4$.  The first part of the lemma
now follows.

For the second part, \Ref{A.10} should be replaced by
\eqs
 && |f(u)g(v)|\pr[L=u]\pr[M=v]\\
&&\quad\left\{ |(1-p-q)^{-(u+v)} - 1| + \Blm\frac{(m-u)\dv}{m\dv} - 1\Brm
     + \Blm \Bl 1 - \frac{pq}{(1-p)(1-q)}\Br^m - 1 \Brm\right\},
\ens
after which  we use the bounds
\[
   \Blm\frac{(m-u)\dv}{m\dv} - 1\Brm \Le \frac{2uv}m;
   \quad \Blm \Bl 1 - \frac{pq}{(1-p)(1-q)}\Br^m - 1 \Brm \Le 4mpq/3.
\]
\ep

\bsk
\begin{lemma}\label{sum-bds}
Let $p_s$, $s\ge j$, be nonnegative numbers summing to $P \le 1$,
and  define
\[
    \s^2_n(r) \ :=\ \ssj (np_s)^r e^{-np_s}, \quad r\ge1; \qquad
    \s^2_n(0) \ :=\ \ssj \min(np_s,1) e^{-np_s} .
\]
Then there exist universal constants $K_r^{(\a)}$, $K_u$, $K_{uv}$ and~$K'$ such that, for any
integers $u\ge v\ge 0$ and for any $\a > 0$,
\eqs
  &(i)&\ssj (np_s)^{u+1} e^{-(1+\a)np_s} \Le K_0^{(\a)} \s^2_n(0); \quad
  (ii)\ \ssj (np_s)^{u+r} e^{-(1+\a)np_s} \Le K_r^{(\a)} \s^2_n(r); \\
  &(iii)&\ssj (np_s)^{u+1} e^{-np_s} \Le K_u nP; \qquad\qquad\
  (iv)\ \Bl \ssj np_s e^{-np_s} \Br^2 \Le K'n\s^2_n(0); \\
  &(v)& \ssj \stj (np_s)^{r+u}(np_t)^{r+v} e^{-n(p_s + p_t)} \Le K_{uv} nP \s^2_n(r).
 % (vi)&& \Bl \ssj (np_s)^{u+r} e^{-np_s} \Br^2 \Le K'n\s^2_n(r)
\ens
\end{lemma}

\proof
The first  inequality reflects the fact that $x^{u+1} e^{-(1+\a)x} \le xe^{-x}$ for
$0\le x\le 1$, whereas $x^{u+1} e^{-(1+\a)x} \le e^{-x} \sup_{z\ge1}\{ze^{-\a z}\}$:
thus we can take $K^{(\a)} = 1/e\a$.  The second is similar in vein, but easier.
The third inequality, and case $u=v=0$ in the fifth, follow from
\[
  \ssj (np_s)^{u+1} e^{-np_s} \Eq n \ssj p_s(np_s)^{u} e^{-np_s} \Le nP (u/e)^u.
\]
For the fifth with $u\ge1$, we write the sum as
\[
  n^2\ssj p_s(np_s)^{r+u-1}e^{-np_s}\,
    \stj p_t[(np_t)^{r+u-1} e^{-np_t}]^{\frac{r+v-1}{r+u-1}}
       \exp\Blb -np_t \frac{u-v}{r+u-1} \Brb,
\]
and use Cauchy--Schwarz to yield the upper bound
\eqs
  \lefteqn{n^2 P \ssj p_s(np_s)^{2r+u+v-2}\exp\Blb -np_s \frac{2r+u+v-2}{r+u-1} \Brb} \\
  &&\Le nP \ssj (np_s)^r e^{-np_s} \max_{x\ge0}\{x^{r+u+v-1}\exp\{-x(r+v-1)/(r+u-1)\}\},
\ens
noting that $r+u-1 \ge 1$.  For the fourth part, Cauchy--Schwarz gives
\[
  \Bl \ssj np_s e^{-np_s} \Br^2 \Le n\ssj np_s e^{-2np_s} \Le \ssj \min\{np_s,e^{-1}\}e^{-np_s}.
\]
% The proof of the last part is much the same.
%\ep

%\bsk
\bigskip
\noindent
{\bf Acknowledgement}
% The authors would like to thank *** for helpful discussions.
This work was carried during a visit to the
Institute for Mathematical Sciences at the National University
of Singapore, whose support is gratefully acknowledged.

\end {document}